\input amstex
\magnification=\magstep1
\documentstyle{amsppt}
\openup1\jot
\hsize6.5truein 
\loadmsbm
\loadeufm
\loadbold
\def\ds{\displaystyle}
\def\re#1{{\text{\rm Re}}(#1)}
\def\gg#1{\Gamma(\displaystyle{#1})}
\vsize9truein
\def\re#1{{\text{\rm Re}}(#1)}
\def\gg#1{\Gamma(\displaystyle{#1})}

\def\eps{\varepsilon}

\def\gen{_m\hskip-1.1pt F_{n}}

\def\hyper{_4\hskip-1.1pt F_3}

\def\gl#1#2{GL(#1,{\Bbb #2})}

\def\Z{{\Bbb Z}}
\def\C{{\Bbb C}}
\def\exp#1{\hbox{e}^{#1}}
\def\R{\Bbb R}
\def\Z{{\Bbb Z}}

\def\eps{\varepsilon}

\def\hyper{_4\hskip-1.1pt F_3}

\def\gl#1#2{GL(#1,{\Bbb #2})}

\def\R{\Bbb R}
\def\Z{{\Bbb Z}}

\parskip8pt\parindent0pt
\topmatter
\title  Coxeter Group Actions on ${_4}F_3(1)$ hypergeometric series 
\endtitle
\author Marc Formichella, R.M. Green and Eric Stade\endauthor
\affil Department of Mathematics \\ University of Colorado \\
Campus Box 395 \\ Boulder, CO  80309-0395 \\ USA \\ {\it  E-mail:}
Marc.Formichella\@Colorado.EDU, rmg\@euclid.colorado.edu, 
Eric.Stade\@colorado.edu \\
\newline
\endaffil

\abstract  We investigate a certain linear combination $K(\vec{x})=K(a;b,c,d;e,f,g)$ of two Saalschutzian hypergeometric series of type ${_4}F_3(1)$.  We first show that $K(a;b,c,d;e,f,g)$ is invariant under the action of a certain matrix group $G_K$, isomorphic to the symmetric group $S_6$, acting on the affine hyperplane $V=\{(a,b,c,d,e,f,g)\in\Bbb C^7\colon e+f+g-a-b-c-d=1\}$. We further develop an algebra of three-term relations for $K(a;b,c,d;e,f,g)$. We show that, for any three elements $\mu_1,\mu_2,\mu_3$  of a certain  matrix group $M_K$, isomorphic to the Coxeter group $W(D_6)$ (of order 23040), and containing the above group $G_K$, there is a relation among $K(\mu_1\vec{x})$, $K(\mu_2\vec{x})$, and $K(\mu_3\vec{x})$, provided no two of the $\mu_j$'s are in the same right coset of $G_K$ in $M_K$.  The coefficients in these three-term relations are seen to be rational combinations of gamma and sine functions in $a,b,c,d,e,f,g$.

The set of $({|M_K|/|G_K|\atop 3})=({32\atop 3})=4960$ resulting three-term relations may further be partitioned into five subsets, according to the Hamming type of the triple $(\mu_1,\mu_2,\mu_3)\ $ in question.  This Hamming type  is defined in terms of Hamming distance between the $\mu_j$'s, which in turn is defined in terms of the expression of the $\mu_j$'s as words in the Coxeter group generators. 

Each three-term relation of a given Hamming type may be transformed into any other of the same type by a change of variable.   An explicit example of each of the five types of three-term relations is provided.   

Our results may be seen to parallel the theory, initiated by Thomae, Whipple, and others, and later put into a group theoretic framework by Beyer, Louck, and Stein, of functional relations among ${_3}F_2(1)$ series.   Moreover, the work of these latter authors and others on transformation formulas for {\it terminating} Saalschutzian ${_4}F_3(1)$ series may be deduced as a special case of our results concerning $K(a;b,c,d;e,f,g)$.    \endabstract

 \endtopmatter

\centerline{\bf Preliminary version, draft 4}

\head  1. Introduction\endhead

The  hypergeometric series$$F(a,b;c;z)=1+\frac{ab}{1!  c}\,z+\frac{a(a+1)b(b+1)}{ 2!  c(c+1)}\,z^2+\frac{a(a+1)(a+2)b(b+1)(b+2)}{3! c(c+1)(c+2)}\,z^3+\cdots\eqno{(1.1)}$$was introduced in 1821 by Gauss [Ga], who developed many of its properties and demonstrated its relationship  to a great variety of elementary and special  functions.  This series, now commonly called a ``Gauss function,'' soon became ubiquitous in mathematics and the physical sciences.

In the latter part of the 1800's, the theory of {\it generalized hypergeometric series} -- series like (1.1), but with arbitrary numbers of ``numerator'' and ``denominator'' parameters (the series (1.1) has two numerator parameters, $a$ and $b$, and one denominator parameter, $c$)   -- began to take form.  From this period through the early part of the 1900's,   properties of -- and, especially, relations among -- these generalized series were studied extensively.  (See [Ba], [Bar1], [Bar2], [T], and [Wh], to name a few.)  

The latter part of the 20th century saw renewed interest  in generalized hypergeometric series---especially those with {\it unit argument}, meaning $z=1$.  This resurgence was due in part to the appearance of  these series in atomic and molecular physics: in the calculation of multiloop Feynman integrals (see [Gr, Chapters 8, 9, and 11]); as $3$-$j$ and $6$-$j$ coefficients in angular momentum theory (cf. [D, Sections 2.7 and 2.9]); and so on.  Generalized hypergeometric series of unit argument have also, recently, commanded a strong presence in the theory of archimedean zeta integrals for automorphic $L$ functions---see, for example, the seminal work of Bump [Bu], and such subsequent works as [ST], [St1], [St2], [St3], and [St4].

Indeed, generalized hypergeometric series have become sufficiently   ubiquitous that it is now typical to drop the adjective ``generalized,'' and we will follow this convention from now on.

Of the relatively recent research in the theory of hypergeometric series, we remark in particular on the article of  [BLS].  In that work, group-theoretic notions were used to explain certain relations among unit argument series  ``of type $_3F_2$,'' meaning series like (1.1), but with an additional parameter in both numerator and denominator (and with $z=1$).  These relations themselves---both ``two-term'' relations and ``three-term'' relations---had been developed much earlier, by Thomae [T], Whipple [Wh], and others.  But the group-theoretic framework for them was new, and it is such a framework that we wish to develop here, in a different but parallel context.

Specifically, we study a certain linear combination $K(a;b,c,d;e,f,g)$  of two unit argument series ``of type $_4F_3$,'' meaning series like (1.1) but with {\it four} numerator parameters and {\it three} denominator parameters (and again with $z=1$). (See Section 2, below, for the precise definition of $K(a;b,c,d,e,f,g)$.)   This particular linear combination of hypergeometric series arises in the theory of archimedean zeta integrals for automorphic $L$ functions.  (See [ST] and [St4].  In the latter work, the function $K(a;b,c,d;e,f,g)$ appears ``disguised'' as a Barnes integral, cf. Proposition 3.1 below.)

The goal of this paper is to describe an ``algebra'' of two-term and three-term relations for $K(a;b,c,d;e,f,g)$, analogous to the $_3F_2$ theory referenced above.  In the case of two-term relations for $K(a;b,c,d;e,f,g)$, we will see that the symmetric group $S_6$ plays a fundamental role.  In the context of the three-term relations, on the other hand, the Coxeter group $W(D_6)$ will appear as the central algebraic object.  Examination of the combinatorial structure of $W(D_6)$ will be crucial to our understanding and categorization of these three-term relations. 

Combinatorial group theory provides a new outlook on the study of relations among hypergeometric series.  We expect that this new perspective will prove fruitful in a variety of contexts similar to those explored in the present work.

The paper [BLS] also provides a group-theoretic framework for certain two-term relations,  originally due to Bailey [Ba, p. 15], for {\it terminating}  $_4F_3$ series of unit argument.  (Such series show up as $6$-$j$ coefficients in angular momentum theory, among other places.)  In the present work (see the latter part of Section 3), we show that these results concerning terminating $_4F_3$ series of unit argument may be derived as special cases of our two-term relations for the series  $K(a;b,c,d;e,f,g)$.

\head  2.  Hypergeometric series: basic definitions and formulas\endhead  

For $a\in\C$, we define
$$(a)_k=(k+a-1)(k+a-2)\cdots
(1+a)a\quad\hbox{$(k\in\Z^+)$;\qquad$(a)_0=1$}.$$Note, by the
functional equation $\gg{s+1}=s\gg{s}$ for the gamma function,
that$$(a)_k={\Gamma(k+a)\over
\Gamma(a)}$$for
$a\ne 0,-1,-2,\ldots$. (Recall that the gamma function has
simple poles at the nonpositive integers.)  Next,  for
$a_1,a_2,\ldots,a_m, e_1,e_2,\ldots,e_{n}, z\in C$, we define
$$\gen(a_1,a_2,\ldots,a_m; e_1,e_2,\ldots,e_{n}; z)=\sum_{k=0}^\infty
{(a_1)_k (a_{2})_k\cdots(a_m)_k\over
k! (e_1)_k (e_2)_k\cdots(e_{n})_k} z^k.$$ (We also sometimes use the notation $${_m}F_n\biggl[{a_1,a_2,\ldots,a_m; \atop e_1,e_2,\ldots,e_{n};} z\biggr]$$for this series.)  If we are to
avoid poles we need to assume that no $e_i$ is a negative
integer or zero; on the other hand if one of the $a_i 's$ is a
nonpositive integer, then the series {\it terminates} (since
then $(a_i)_k$ will be zero for $k> -a_i$).  

The above series is
called a {\it hypergeometric series of type $\gen$}; we will restrict our attention
to the case
$n=m-1$.  The series $_mF_{m-1}$ converges absolutely for $|z|<1$, or for
$|z|=1$ and $\re{\sum f_i-\sum e_i}>0$. (See [Ba, Chapter 2].) If $z=1$, this series is said to be ``of
type
$_mF_{m-1}(1)$'' or ``of unit argument.'' If $\re{\sum f_i-\sum e_i}=1$, the
series is said to be ``Saalschutzian.''

Our primary object of interest, in this paper, will be the following linear combination $K(a;b,c,d;e,f,g)$ of Saalschutzian $\hyper(1)$ series:
$$\displaylines{K(a;b,c,d;e,f,g)=\frac{{_4}F_3 (a,b,c,d;e,f,g;1)}{\Gamma(e)\Gamma(f)\Gamma(g)\gg{1+a-e}\gg{1+a-f}\gg{1+a-g}} \cr
 +\frac{{_4}F_3 (a,1+a-e,1+a-f,1+a-g;1+a-b,1+a-c,1+a-d;1)}{\gg{1+a-b}\gg{1+a-c}\gg{1+a-d}\gg{b}\gg{c}\gg{d}}\cr=\frac{1}{\gg{a}\gg{b}\gg{c}\gg{d}\gg{1+a-e}\gg{1+a-f}\gg{1+a-g}}\cr
\cdot\bigl[{_4}F_3^* (a,b,c,d;e,f,g;1) +{_4}F_3^* (a,1+a-e,1+a-f,1+a-g;1+a-b,1+a-c,1+a-d;1)],\cr}$$where $g=1+a+b+c+d-e-f$ and, by definition,  $$\displaylines{
{_4}F_3^* (a,b,c,d;e,f,g;1)={\Gamma(a)\Gamma(b)\Gamma(c)\Gamma(d)
\over\Gamma(e)\Gamma(f)\Gamma(g)}
{\hyper(a,b,c,d;e,f,g;1)}\cr=\sum_{k=0}^\infty
{\Gamma(k+a)\Gamma(k+b)\Gamma(k+c)\Gamma(k+d)\over
k!\Gamma(k+e)\Gamma(k+f)\Gamma(k+g)}   .}$$We will also sometimes denote $K(a;b,c,d;e,f,g)$ by $\displaystyle{K\left[{a;b,c,d;\atop e,f,g}\right]}$.

Again, the gamma function never vanishes; it follows that $K(a;b,c,d;e,f,g)$ is entire in each of its independent variables.

We will prove the following concerning relations among  $K(a;b,c,d;e,f,g)$ series.

(i)  {There is a group $G_K$ of linear transformations of the affine hyperplane $$V=\{\vec{x}=(a,b,c,d,e,f,g)^T\in\C^7\,\colon\, e+f+g-a-b-c-d=1\},\eqno{(2.1)}$$such that $G_K$ is isomorphic to the symmetric group $S_6$;
$G_K$ fixes the first coordinate of $V$ ($\gamma(a,b,c,d,e,f,g)^T=(a,\gamma(b),\gamma(c),\gamma(d),\gamma(e),\gamma(f),\gamma(g))^T$ for all $\gamma\in G_K$); 
and $$
  K(g\vec{x})
=K(\vec{x})
$$for all $g\in G_K$.}   

 {(ii)}  {There is a group $M_K$ of linear transformations of $V$, such that $M_K$ is isomorphic to the Coxeter group $W(D_6)$ (or order 23040) in such a way that
the subgroup $G_K$ of $M_K$ can be identified with the parabolic subgroup
$W(A_5)$.  Furthermore, for any $\mu_1,\mu_2,\mu_3\in M_K$ such that no two of these elements are in the same right coset of $M_K$ modulo $G_K$,
there exists an expression of the form$$ \gamma_1(\vec{x}) K(\mu_1{\vec x})+\gamma_2(\vec{x}) K(\mu_2 {\vec x})+\gamma_3(\vec{x})  K(\mu_3 {\vec x}) =0,\eqno{(2.2)}$$where  $\gamma_1$, $\gamma_2$, $\gamma_3$ are ratios of gamma functions in $a,b,c,d,e,f,g$.
(Because our coset space $G_K\backslash M_K$  has order
$|M_K|/|G_K|=23040/720=32$, there are $({32\atop 3})=4960$ ``different'' three-term relations of the kind just described.)} 

{(iii)}  {The induced action of $M_K$ on distinct triples of right
$G_K$-cosets has precisely five orbits.
The three-term relations referenced above fall into five corresponding
families, according to the ``Hamming type'' of the triple $(\mu_1,\mu_2,\mu_3)$.  This Hamming type is defined in terms of a metric, which we call the ``Hamming distance,'' on the coset space $G_K\backslash M_K$. 
 
{(iv)}   {In a relation (2.2), suppose the Hamming distance between two of the coset representatives $\mu_j$ and $\mu_k$ ($ j,k\in\{1,2,3\}$ and $j\ne k$) is $2n$:  then the ``opposite'' coefficient $\gamma_ \ell$ ($ \ell\in\{1,2,3\}-\{j,k\}$) may be written as a sum  of $2^{n-1}$ ratios of gamma and sine  functions in $(a,b,c,d,e,f,g)$.}

{(v)}  {Any three-term relation of a given Hamming type may be transformed into any other  of the {\it same} type by a change of variable $$(a,b,c,d,e,f,g)\to \rho(a,b,c,d,e,f,g)\quad(\rho\in M_K)$$ applied to all elements (including the coefficients) of the first relation.}
}

Our investigations below  will require some classical facts concerning the Gauss function
$F(a,b;c;z)$.  It is well-known that this function has continuation in
$z$, and satisfies the functional equations:
$$F(a,b;c;z)=(1-z)^{-b}
 F\biggl(c-a,b;c;\frac{z}{z-1}\biggr);\eqno{(2.3)}$$
$$\aligned 
F(a,b;c;z)=&{\Gamma(c) \Gamma(c-a
-b)\over
\Gamma(c-a) \Gamma(c-b)} 
F(a,b;1+a+b-c;1-z)  \\+ (1-z)^{c-a-b}&
{\Gamma(c) \Gamma(a+b-c)\over
\Gamma(a) \Gamma(b)} 
F(c-a,c-b;1+c-a-b;1-z)\endaligned\eqno{(2.4)}$$(among
others). 

We will also make use of the notion of a {\it Mellin--Barnes  integral}, by which we mean one of the
form$$\int_t\prod_{k=1}^m
\Gamma^{\eps_k}(a_k+t)\prod_{\ell=1}^n
\Gamma^{\eps_\ell}(b_\ell-t) \,z^t \,dt, $$where $m,n$
are nonnegative integers; ${\eps_k},{\eps_\ell}=\pm1$;
$a_k,b_\ell,z\in{\Bbb C}$.   (Such an integral is, in certain
contexts,  also called a {\it Meijer's $G$-function}.) The path of
integration is a line parallel to the imaginary axis, indented
if necessary to ensure that any poles of $\prod_{k=1}^m
\Gamma^{\eps_k}(a_k+t)$ are to the left of this path, while 
any
poles of $\prod_{\ell=1}^n
\Gamma^{\eps_\ell}(b_\ell-t)$ are to its right. 
(Note that poles arise only for those $k$ with $\eps_k=1$,
and those $\ell$ with $\eps_\ell=1$, since the gamma function is
never zero.  Note also the tacit assumption that, if
$\eps_k=\eps_\ell=1$, then
$a_k+b_\ell\not\in{\Bbb Z}$; otherwise the path of integration
could not be chosen as above.)  In this paper,  {\it the path of integration of any 
Mellin--Barnes  integral will always be of this form}.

By moving the above line of integration to the left or
to the right and summing residues, using the facts
that$$\hbox{Res}_{s=-n}\gg{s}={(-1)^n\over n!}\quad(n=0,1,2,\ldots)\eqno{(2.5)}$$and 
$$\gg{s}={\pi\over
\gg{1-s}\sin\pi s}\quad (s\not\in{\Bbb Z}),\eqno{(2.6)}$$ one may at least formally
evaluate a  Mellin--Barnes integral as a linear combination of
powers of $z$ times hypergeometric series (in $\pm z^{\pm1}$). 
Of course doing so requires absolute convergence of the relevant
integrals and series, and that integrals over certain paths
approach zero as the paths move out to infinity, etc. 

This connection between Mellin--Barnes integrals and hypergeometric series will be exploited at various points throughout this paper.   In particular, the following lemma will be of use in Section 3 below.

\proclaim{Lemma 2.1}  For $|\hbox{\rm arg}(z)|<\pi$, we have$$\displaylines{\frac{1}{2\pi i}\int_t \gg{a+t}\gg{b+t}\gg{c-t}\gg{d-t}\,z^t\,dt\cr=\cases \displaystyle z^{c}\frac{\gg{a+c}\gg{b+c}\gg{a+d}\gg{b+d}}{\gg{a+b+c+d}} F\biggl[{a+c,b+c;\atop a+b+c+d;}1-z\biggr] & \text{if $|z|<1$},\\\displaystyle
z^{-a}\frac{\gg{a+c}\gg{b+c}\gg{a+d}\gg{b+d}}{\gg{a+b+c+d}} F\biggl[{a+c,a+d;\atop a+b+c+d;}1-z^{-1}\biggr]&\text{if $|z|>1$.}\\\endcases}$$
\endproclaim

\demo{Proof}  In the case $|\hbox{arg}(z)|<\pi$ and $|z|<1$, we may move the line of integration all the way to the right, and sum residues at the poles of the integrand.  The validity of such an approach follows from arguments like those in [Ba, Section 1.5], and makes use of {\it Stirling's formula}
$$\gg{s}\sim
\sqrt{2\pi} s^{s-1/2}\exp{-s}, $$which holds
uniformly as
$|s|\rightarrow\infty$ in any sector of the form
$|\hbox{arg}(s)| <\pi-\delta$. 

 By our tacit assumptions  concerning the positioning of the path of integration, and by (2.5) and (2.6), we have$$\displaylines{\frac{1}{2\pi i}\int_t \gg{a+t}\gg{b+t}\gg{c-t}\gg{d-t}\,z^t\,dt\cr=\sum_{n=0}^\infty\frac{(-1)^n\gg{n+a+c}\gg{n+b+c}\gg{d-c-n}}{n! }\,z^{n+c}
\cr+\sum_{n=0}^\infty\frac{(-1)^n\gg{n+a+d}\gg{n+b+d}\gg{c-d-n}}{n! }\,z^{n+d}
\cr=\frac{\pi z^{c}}{\sin\pi (d-c)}\sum_{n=0}^\infty\frac{ \gg{n+a+c}\gg{n+b+c} }{n!\gg{n+1+c-d}}\,z^{n}
\cr +\frac{\pi z^{d}}{\sin\pi (c-d)}\sum_{n=0}^\infty\frac{ \gg{n+a+d}\gg{n+b+d} }{n!\gg{n+1+d-c}}\,z^{n}
\cr=  z^{c}\gg{a+c}\gg{b+c} \gg{d-c}F(a+c,b+c;1+c-d;z) \cr
 + z^{d}\gg{a+d}\gg{b+d} \gg{c-d}F(a+d,b+d;1+c-d;z) \cr=z^{c}\frac{\gg{a+c}\gg{b+c}\gg{a+d}\gg{b+d}}{\gg{a+b+c+d}} F(a+c,b+c;a+b+c+d;1-z), }$$the last step by  (2.4). So we are done in the case $|z|<1$.

For the case $|z|>1$, we simply make the change of variable $t\to -t$ in the integral in question, and proceed essentially as above.\qed  

We remark that, under appropriate conditions on $a,b,c,d$, the above arguments extend to the case $z=1$.  Since $F(\alpha,\beta;\gamma;0)=1$, we conclude that
$$\frac{1}{2\pi i}\int_t \gg{a+t}\gg{b+t}\gg{c-t}\gg{d-t}\,dt=\frac{\gg{a+c}\gg{b+c}\gg{a+d}\gg{b+d}}{\gg{a+b+c+d}} ,$$under these conditions.  This result is known as {\it Barnes' Lemma}, and is shown to hold as long as none of the quantities $a+c,b+c,a+d,b+d$ is an integer (so that sequences of poles in the integrand may be separated in the usual way).
 \enddemo
In the case $z=1$, a Mellin--Barnes  integral is sometimes called, simply, a {\it Barnes integral}.   In the next section, we will provide a useful Barnes integral representation of our function $K(a;b,c,d;e,f,g)$.

\head  3.  Two-term relations\endhead

In this section, we explore ``two-term relations'' for $K(a;b,c,d;e,f,g)$, meaning invariances of this function under certain linear transformations on the affine hyperplane $V$ defined in (2.1).  Central to the study of such invariances will be the following integral representation of $K(a;b,c,d;e,f,g)$.
\noindent \proclaim{Proposition 3.1} 
If $e+f+g-a-b-c-d=1$, 
then 
$$\displaylines{  K\biggl[{a;b,c,d;
\atop
e,f,g} \biggr]
\cr=\frac{1}{ {
\gg{1+a-e}\gg{1+a-f}\gg{b}\gg{c}\gg{e-b}\gg{e-c}\gg{f-b}\gg{f-c}
\gg{g-d} }}
\cr\cdot{1\over 2\pi i}
\int_t \frac{\biggl[{\displaystyle\gg{g-d+t}\gg{1+a-e+t}\gg{1+a-f+t}
\atop \displaystyle\cdot\gg{d+b-g-t}
\gg{d+c-g-t}\gg{-t}}\biggr]}{\gg{d-t}\gg{1+a-d+t}} \,dt  .}$$
\endproclaim

\demo{Proof} This result is essentially Proposition 1 of [ST], with a change of variable. For completeness, we provide a similar though somewhat simplified proof here, but leave out technical details concerning convergence, analytic continuation, and so on.

 We denote by
$I(a,b,c,d,e,f,g)$ the integral on the right-hand side of our proposition (including the factor of $(2 \pi i)^{-1}$). Using the classical integral formula $$\frac{\gg{\alpha}\gg{\beta}}{\gg{\alpha+\beta}}=\int_0^\infty x^\alpha (1+x)^{-(\alpha+\beta)}\,\frac{dx}{x}\quad (\hbox{Re}(\alpha),\hbox{Re}(\beta)>0),\eqno{(3.1)}$$we may write 
$$\displaylines{  I(a,b,c,d,e,f,g)={1\over\gg{1+a-g}\gg{d}}
 \cr \cdot{1\over
2\pi i}\int 
\gg{1+a-e+t}\gg{1+a-f+t} \gg{d+b-g-t} \gg{d+c-g-t}
\cr\cdot
\int_0^\infty\hskip-5pt\int_0^\infty x^{1+a-g} (1+x)^{-(1+a-d+t)} 
y^{d} (1+y)^{-(d-t)}\,
{dx\over x}{dy\over y}\,dt\cr ={1\over\gg{1+a-g}\gg{d}}\int_0^\infty\hskip-5pt\int_0^\infty x^{1+a-g} (1+x)^{-(1+a-d)} 
y^{d} (1+y)^{-d}
 \cr\cdot{1\over
2\pi i}\int 
\gg{1+a-e+t}\gg{1+a-f+t} \gg{d+b-g-t} \gg{d+c-g-t}\biggl(\frac{1+y}{1+x}\biggr)^t\hskip-2pt
 dt {dx\over x}{dy\over y}.}$$ 

Now in the above expression for  $I(a,b,c,d,e,f,g)$, we break the integral in $x$ and $y$ into
an integral over $0<y<x<\infty$, plus an integral over 
$0<x<y<\infty$. Each of the resulting integrals may be evaluated using Lemma 2.1; the upshot is that  
$$\displaylines{ I(a,b,c,d;e,f,g)=\frac{\gg{e-b}\gg{e-c}\gg{f-b}\gg{f-c}}{\gg{1+a-g}\gg{d}\gg{1+a+d-g}}
 \cr \cdot
\biggl\lbrace\int_0^\infty\hskip-5pt\int_y^\infty x^{1+a-g} (1+x)^{-(1+a+b-g)} 
y^{d} (1+y)^{-(g-b)}  F\biggl(f-c,e-c;1+a+d-g;\frac{x-y}{1+x}\biggr) 
{dx\over x}{dy\over y} 
\cr+\int_0^\infty\hskip-5pt\int_x^\infty x^{1+a-g} (1+x)^{-(e-d)} 
y^{d} (1+y)^{-(1+a+d-e)}  F\biggl(f-c,f-b;1+a+d-g;\frac{y-x}{1+y}\biggr) 
{dy\over y} {dx\over x} \biggr\rbrace.}$$ 

Into the first of the integrals on the right, we now  substitute
$x\to \frac{u+y}{1-u} $, and follow this with the substitution $y\to y u$;   into the second, we substitute 
$y\to \frac{u+x}{1-u}$, and follow this with the substitution $x\to x u$.
  The result,
after some rearranging,
is
that$$\displaylines{ I(a,b,c,d;e,f,g)=\frac{\gg{e-b}\gg{e-c}\gg{f-b}\gg{f-c}}{\gg{1+a-g}\gg{d}\gg{1+a+d-g}}
 \cr \cdot
\biggl\lbrace\int_0^\infty  \hskip-2pt y^{d} (1+y)^{-(g-a)} \hskip-4pt \int_0^1 
\hskip-2pt u^{a+d-g} (1-u)^{b-1}  (1+yu)^{-a} 
 F (f-c,e-c;1+a+d-g;u )\,
du {dy\over y} 
\cr+\int_0^\infty \hskip-2pt x^{1+a-g} (1+x)^{-(1-d)}\hskip-4pt \int_0^1 \hskip-2pt
u^{a+d-g} (1-u)^{a-e} (1+xu)^{-a} F(f-c,f-b;1+a+d-g;u)\,
du {dx\over x} \biggr\rbrace.}$$To the integrals in $u$ on the right, we may apply the formula$$\displaylines{\int_0^1 u^{\gamma-1} (1-u)^{\rho-1}(1-zu)^{-\sigma}F(\alpha,\beta;\gamma;u)\,du\cr=\frac{\gg{\gamma}\gg{\rho}\gg{\gamma+\rho-\alpha-\beta}}{\gg{\gamma+\rho-\alpha}\gg{\gamma+\rho-\beta}}(1-z)^{-\sigma}\,{_3}F_2\biggl(\rho,\sigma,\gamma+\rho-\alpha-\beta;\gamma+\rho-\alpha,\gamma+\rho-\beta;\frac{z}{z-1}\biggr)\cr(\hbox{Re}(\gamma),\hbox{Re}(\rho),\hbox{Re}(\gamma+\rho-\alpha-\beta)>0;\ |\hbox{arg}(1-z)|<\pi),}$$cf. equation 7.512(9) in [GR]. The result is that 
$$\displaylines{ I(a,b,c,d;e,f,g)=\frac{\gg{e-b}\gg{e-c}\gg{f-b}\gg{f-c}}{\gg{1+a-g}\gg{d} }
 \cr \cdot
\biggl\lbrace\frac{\gg{b}\gg{c}}{\gg{e}\gg{f}}\int_0^\infty y^{d} (1+y)^{-g} \,{_3}F_2\biggl(b,a,c;e,f;\frac{y}{1+y}\biggr){dy\over y}+ \frac{\gg{1+a-e}\gg{1+a-f}}{\gg{1+a-b}\gg{1+a-c}}
\cr\cdot\int_0^\infty x^{1+a-g} (1+x)^{-(1+a-d)}\,{_3}F_2\biggl(1+a-e,a,1+a-f;1+a-b,1+a-c;\frac{x}{1+x}\biggr){dx\over x} \biggr\rbrace.}$$Expanding out   the above ${_3}F_2$'s as infinite series, exchanging the summations with the integrations, and performing the remaining integrals using (3.1), we get
$$\displaylines{ I(a,b,c,d;e,f,g)
 \cr =\frac{\gg{e-b}\gg{e-c}\gg{f-b}\gg{f-c}\gg{g-d}}{\gg{a}\gg{1+a-g} \gg{d}}
\biggl\lbrace \sum_{k=0}^\infty
{\Gamma(k+a)\Gamma(k+b)\Gamma(k+c)\Gamma(k+d)\over
k!\Gamma(k+e)\Gamma(k+f)\Gamma(k+g)} \cr + \sum_{k=0}^\infty{\Gamma(k+a)\Gamma(k+1+a-e)\Gamma(k+1+a-f)\Gamma(k+1+a-g)\over
k!\Gamma(k+1+a-b)\Gamma(k+1+a-c)\Gamma(k+1+a-d)} \biggr\rbrace,}$$which amounts precisely to the statement of the lemma, and we are done.
\qed
\enddemo 

From the above proposition we may deduce the following fundamental two-term relation.

\proclaim{Corollary 3.2} If $e+f+g-a-b-c-d=1$, then 
$$\displaylines{K\biggl[{a;b,c,d;\atop
e,f,g}\hskip-1pt \biggr] 
=K\biggl[{a;e-c,e-b,d;\atop
 e,1+a+d-g,1+a+d-f}  \biggr ].}$$\endproclaim
\demo{Proof} By Proposition 3.1, we have$$\displaylines{  K\biggl[{a;e-c,e-b,d;
\atop
e,1+a+d-g,1+a+d-f} \biggr]
\cr=\frac{1}{ \biggl[{\ds
\gg{1+a-e}\gg{g-d}\gg{e-c}\gg{e-b}\gg{c}\gg{b}\atop\ds\cdot\gg{1+a+c+d-e-g}\gg{1+a+b+d-e-g}
\gg{1+a-f} }\biggr]}
\cr\cdot{1\over 2\pi i}
\int \frac{\biggl[{\ds\gg{1+a-f+t}\gg{1+a-e+t}\gg{g-d+t} \atop \ds\cdot\gg{e+f-a-c-1-t}
\gg{e+f-a-b-1-t}\gg{-t}}\biggr]}{\gg{d-t}\gg{1+a-d+t}} \,dt
\cr=\frac{1}{ {
\gg{1+a-e}\gg{g-d}\gg{e-c}\gg{e-b}\gg{c}\gg{b}\gg{f-b}\gg{f-c}
\gg{1+a-f} }}
\cr\cdot{1\over 2\pi i} 
\int {\biggl[{\ds\gg{1+a-f+t}\gg{1+a-e+t}\gg{g-d+t}\atop\ds\cdot\gg{d+b-g-t}
\gg{d+c-g-t}\gg{-t}}\biggr]\over\gg{d-t}\gg{1+a-d+t}} \,dt  ,}$$the last step by the hypothesized   relation among $a,b,c,d,e,f,g$.  But by Proposition 3.1 again, the right-hand side above equals $K(a;b,c,d;e,f,g)$, and we are done.\qed
\enddemo

The function $K(a;b,c,d;e,f,g)$ is also, by its very definition, invariant under permutations of $b,c,d$ and of $e,f,g$.  The group of transformations generated by these permutations and by the change of variable
$$\biggl({a,b,c,d,\atop e,f,g}\biggr)\to \biggl({a,e-c,e-b,d,
\atop
e,1+a+d-g,1+a+d-f}\biggr)$$tacit in the above corollary is described in the  following proposition.

\proclaim{Proposition 3.3}
   Let $Q_6$ denote the group of $6\times 6$ permutation matrices (matrices with exactly one entry of 1 in each row and column, and zeroes elsewhere). Denote by $P_{1,6}$ the group of all matrices $p\in SL(7,\Bbb R)$ of the form$$p=\pmatrix 1& 0 \\
0 &q  \endpmatrix,$$where $q\in Q_6$, so that $P_{1,6}$ is isomorphic to the symmetric group $S_6$. 

 For complex numbers $a,b,c,d,e,f,g$ satisfying $e+f+g-a-b-c-d=1$,  let $\vec{x}$ denote  the column vector $(a,b,c,d,e,f,g)^T,$ and write $K(\vec{x})$ for $K(a;b,c,d;e,f,g).  $   Also put
$$S=\pmatrix
  1 & 0 & 0 & 0 & 0 & 0 & 0  \\
 0 & 0 & 1 & 1 & 0 & 0 & 0 \\
 0 & 1 & 0 & 1 & 0 & 0 & 0 \\
 0 & 1 & 1 & 0 & 0 & 0 & 0 \\
 0 & 1 & 1 & 1 & 1 & 0 & 0  \\
 0 & 1 & 1 & 1 & 0 & 1 & 0  \\
 0 & 1 & 1 & 1 & 0 & 0 & 1   \endpmatrix\in GL(7,\Bbb R).$$Then $$K(g\vec{x})=K(\vec{x})$$ for all $g$ in the group $G_K=\{ S pS^{-1}\colon p\in P_{1,6}\}$.  (Here, $g\vec{x}$ denotes matrix multiplication.)
\endproclaim

\demo{Proof}   As observed above, $K(\vec{x})$ is invariant under any element of $P_{1,6}$ that permutes the second through the fourth coordinates of $\vec{x}$ among themselves, and/or does the same to the last three coordinates.  In particular, it is invariant under the ``transposition'' matrices $(23)$, $(34)$, $(56)$, and $(67)$.  (Here and throughout, we identify an element of $P_{1,6}$, or more generally any element of the group $Q_7$  of $7\times  7$ permutation matrices, with the permutation it induces on the coordinates of $\vec{x}$.  For example, $(156)\vec{x}=(f,b,c,d,a,e,g)$.)

Direct computation shows that, for $S$ as above, $S(23)S^{-1}=(23)$, $S(34)S^{-1}=(34)$, $S(56)S^{-1}=(56)$, and $S(67)S^{-1}=(67)$.  But also by direct computation, we have$$S(45)S^{-1}\vec{x}= \biggl({a,e-c,e-b,d,\atop
 e,1+a+d-g,1+a+d-f}\biggr),$$so that, by the above corollary, $K(\vec{x})$ is also invariant under $S(45)S^{-1}$.  Then $K(\vec{x})$ is in fact invariant under $\{S(i,i+1)S^{-1}\colon 2\le i\le 6 \}$, and consequently under the group generated by this set.  This group equals $G_K$, since the transposition matrices $(23)$, $(34)$, $(45)$, $(56)$, and $(67)$ generate $P_{1,6}$, and we are done. \qed
\enddemo

\remark{Remark 3.4}  Note that the first coordinate of $g\vec{x}$ equals that of $\vec{x}$, for all $g\in G_K$.\endremark

It should also be noted that our $K(a;b,c,d;e,f,g)$ relations above reduce to known relations concerning {\it terminating} $_4F_3(1)$ hypergeometric series, in certain special cases.
For example, let us start with Corollary 3.2, written out explicitly in terms of ${_4}F_3(1)$ series:
$$\displaylines{\frac{{_4}F_3 (a,b,c,d;e,f,g;1)}{\Gamma(e)\Gamma(f)\Gamma(g)\gg{1+a-e}\gg{1+a-f}\gg{1+a-g}} \cr
 +\frac{{_4}F_3 (a,1+a-e,1+a-f,1+a-g;1+a-b,1+a-c,1+a-d;1)}{\gg{1+a-b}\gg{1+a-c}\gg{1+a-d}\gg{b}\gg{c}\gg{d}}\cr=
\frac{{_4}F_3({a,e-c,e-b,d; 
 e,1+a+d-g,1+a+d-f} )}{\gg{e}\gg{1+a+d-g}\gg{1+a+d-f}\gg{1+a-e}\gg{g-d}\gg{f-d}}\cr +\frac{{_4F_3}({a,1+a-e,g-d,f-d;1+a+c-e,1+a+b-e,1+a-d})}{
\gg{1+a+c-e}\gg{1+a+b-e}\gg{1+a-d}\gg{e-c}\gg{e-b}\gg{d}}.}$$Now let us allow $d$ to approach $ -n$, where $n$ is a nonnegative integer.    Since$$\lim_{d\to -n} \frac{\Gamma(k+d)}{\Gamma(d)}
=(-n)_k\hbox{ and } \lim_{d\to -n} \frac{1}{\Gamma(d)}=0$$for nonnegative integers $k$ and $n$, we get
$$\aligned{_4}F_3\biggl[{a,b,c,-n;\atop e,f,g;}1\biggr] =&
\frac{\Gamma(f)\Gamma(g)\Gamma(1+a-f) \Gamma(1+a-g)}{\Gamma(1+a-f-n)\Gamma(1+a-g-n) \Gamma(f+n)\Gamma(g+n)}\\ 
\cdot &{_4F_3}\biggl[{a,e-c,e-b,-n;\atop e,1+a-g-n,1+a-f-n ;}1\biggr] .\endaligned$$This latter relation is just equation (2.3) in [BLS] (with some of the numerator and denominator parameters on either side permuted).  Note that each of the $_4F_3(1)$ series here  terminates, since  again $(-n)_k=0$ for $k$ a nonnegative integer greater than $n$.

The above identity implies that the Saalschutzian series  
$$\displaylines{ (e)_n(f)_n (g)_n \,_4F_3\biggl[{a,b,c,-n;\atop e,f,g;}1\biggr] }$$is invariant under the change of variables $$(a,b,c,e,f,g)\to (a,e-c,e-b,e,1+a-g-n,1+a-f-n).\eqno{(3.2)}$$ We also have trivial invariance under permutations of $(a,b,c)$ and $(e,f,g)$.  It is shown in [Beyer] that the group of transformations generated by (3.2) and by these permutations is isomorphic to $S_6$.   What we have shown here is that this group of symmetries of the series $ (e)_n(f)_n (g)_n \,_4F_3\bigl[{a,b,c,-n;\atop e,f,g;}\scriptstyle1\bigr]$ follows readily from our two-term relations for the functions $K(a;b,c,d;e,f,g)$ above.

\head   4.   Coxeter groups and linear transformations\endhead

\definition{Definition 4.1}  
We denote by $M_K$ the subgroup of $SL(7,\R)$ generated by $G_K$ (cf. Proposition 3.3 above) and by the transposition matrix (12).
\enddefinition

One of our ultimate goals in this work is to derive, for any $\mu_1,\mu_2,\mu_3$ in the above group $M_K$ such that no two $\mu_j$'s are in the same right coset of $M_K$ modulo $G_K$, a relation of the form$$\gamma_1(\vec{x}) K(\mu_1\vec{x})+\gamma_2(\vec{x}) K(\mu_2\vec{x})
+\gamma_3(\vec{x}) K(\mu_3\vec{x})=0,$$where the $\gamma_j(\vec{x})$'s are certain  rational combinations of gamma functions and sine functions (whose arguments are linear combinations of the coordinates of $\vec{x}$). To develop, and to describe in detail, the ``algebra'' of such three-term relations, we will first need to investigate the combinatorial structure of $M_K$.  In this section and in the next, we present a sequence of results to this end.

\definition{Definition 4.2}
The {\it Dynkin diagram $G(D_6)$ of type $D_6$} is a graph with vertex set 
$V(D_6) = \{1', 1, 2, 3, 4, 5\}$.  Two vertices $i$ and $j$ in the subset 
$\{1, 2, 3, 4, 5\}$ are adjacent (connected by an edge) if and only if 
$|i - j| = 1$, and the vertex $1'$ is adjacent only to $2$.  

If $i$ and $j$ are any two vertices of $V(D_6)$, we define the integer $$
m(i, j) = \cases 
1 & \text{ if $i = j$,}\cr
2 & \text{ if $i$ is not adjacent to $j$, and}\cr
3 & \text{ if $i$ is adjacent to $j$.}\cr
\endcases $$

The {\it Coxeter group $W = W(D_6)$ of type $D_6$} is the group given
by the presentation $$
\langle s_{1'}, s_1, s_2, s_3, s_4, s_5 \ | \ (s_i s_j)^{m(i, j)} = 1 \rangle
.$$
\enddefinition

\definition{Definition 4.3}
Let $1 \leq i < j \leq N$ and regard $(i, j)$ as a transposition in the
symmetric group $S_N$.  For each $\varepsilon \in \{\pm 1\}$, we define 
$M_N^\varepsilon(i, j)$ to be the matrix
$A = (a_{ij}) \in M_N({\Bbb R})$ given by $$
a_{kl} = \cases 
1 & \text{ if $l = (i, j).k$ and $k = l$;}\cr
\varepsilon & \text{ if $l = (i, j).k$ and $k \ne l$;}\cr
0 & \text{ otherwise.}\cr
\endcases
$$
\enddefinition

\proclaim{Lemma 4.4}
\item{}{\rm (i)} {There is a unique homomorphism of groups $$
\phi : W(D_6) \rightarrow GL_6({\Bbb R})
$$ such that $\phi(s_i) = M^+_6(i, i+1)$ for $i \in \{1, 2, 3, 4, 5\}$,
and $\phi(s_{1'}) = M^-_6(1, 2)$.}
\item{}{\rm (ii)} {The map $\phi$ is injective, and
the image of $\phi$ consists of all monomial matrices with entries in
$\{+1, -1, 0\}$ that have an even number of occurrences of $-1$.}
\item{}{\rm (iii)} {The group $W(D_6)$ acts faithfully and transitively on 
the set $\Omega$ of $32$ vectors of the form $$
(\pm 1, \pm 1, \pm 1, \pm 1, \pm 1, \pm 1)^T
$$ that contain an even number of occurrences of $-1$.}
\endproclaim

\demo{Proof}
This is a restatement of the standard semidirect product construction of the
Coxeter groups of type $D_n$ in the special case $n = 6$; see \cite{H, 
\S2.10} for more details.
\qed\enddemo

\proclaim{Lemma 4.5}
\item{}{\rm (i)} {The order of the group $W(D_6)$ is $2^{6-1} . 6! = 23040$.}
\item{}{\rm (ii)} {The stabilizer in $W(D_6)$ of the element $$
\omega_0 = (1, 1, 1, 1, 1, 1)^T \in \Omega
$$ is precisely the subgroup generated by
the set $S' = \{s_1, s_2, s_3, s_4, s_5\}$.}
\item{}{\rm (iii)} {The stabilizer in $W(D_6)$ of any element $\omega \in \Omega$
is isomorphic to the symmetric group $S_6$.}
\endproclaim

\demo{Proof}
Part (i) is well known (see \cite{H, \S2.11}) and follows, for example,
from enumerating the matrices in Lemma 4.4.

To prove (ii), we first note that the generators listed all fix $\omega_0$.
The subgroup $W_{S'} \leq W$ generated by $S'$ is isomorphic to $S_6$, 
because it is a parabolic subgroup of $W(D_6)$ of type $A_5$ 
(see \cite{H, \S1.10}).  The index of $W_{S'}$ in $W$ is $23040/6! = 32$,
and the fact that $W(D_6)$ acts transitively on $\Omega$ (Lemma 4.4 (iii))
then shows that $W_{S'}$ is the full stabilizer in $W$ of $\omega$, by the
orbit-stabilizer theorem.

Since $W(D_6)$ acts transitively on $\Omega$, the stabilizers of elements
of $\Omega$ must be conjugate subgroups of $W(D_6)$.  Part (iii) now follows 
from part (ii).
\qed\enddemo

\proclaim{Lemma 4.6}
The subgroup $G_1$ of $GL_7({\Bbb R})$ that is generated by the matrices $$
M^+_7(2, 3), 
M^+_7(3, 4), 
M^+_7(4, 5), 
M^+_7(5, 6), 
M^+_7(6, 7)
$$ and $$
A_1 = \left( \matrix
1 & 1 & 1 & 0 & 0 & 0 & 0\cr
0 & 0 & -1 & 0 & 0 & 0 & 0\cr
0 & -1 & 0 & 0 & 0 & 0 & 0\cr
0 & 0 & 0 & 1 & 0 & 0 & 0\cr
0 & 0 & 0 & 0 & 1 & 0 & 0\cr
0 & 0 & 0 & 0 & 0 & 1 & 0\cr
0 & 0 & 0 & 0 & 0 & 0 & 1\cr
\endmatrix \right)
$$ is isomorphic to $W(D_6)$.  Furthermore, the matrices can be respectively
identified with the generators $$
s_1, s_2, s_3, s_4, s_5, s_{1'}
$$ of $W(D_6)$.
\endproclaim

\demo{Proof}
We first observe that the first column of every generator of $G_1$ is
the vector $(1, 0, 0, 0, 0, 0, 0)^T$, and it follows from this that every
element of the group $G_1$ has the same property.  Because of this, 
there exists a homomorphism of groups $\psi_1 : G_1 \rightarrow GL_6({\Bbb R})$, 
where $\psi_1(g)$ is defined to be the matrix obtained by removing the 
first row and column of $g$.  The images of the generators of $G_1$ under
$\psi_1$ are precisely the matrices $$
M^+(1, 2), 
M^+(2, 3), 
M^+(3, 4), 
M^+(4, 5), 
M^+(5, 6), 
M^-(1, 2)
$$ in $GL_6({\Bbb R})$.  The latter group of matrices is isomorphic to $W(D_6)$
by Lemma 4.4 (i), so the map $\phi^{-1} \circ \psi_1$ is a homomorphism of
groups from $G_1$ to $W(D_6)$ that identifies the generators of $G_1$ with
$s_1, s_2, s_3, s_4, s_5$ and $s_{1'}$, respectively.

To complete the proof, it suffices to show that $\phi^{-1} \circ \psi_1$ is
an isomorphism, and we do this by constructing its inverse, namely a 
group homomorphism $\psi_2 : W(D_6) \rightarrow G_1$ such that $$
\psi_2(s_i) = \cases
M^+_7(i+1, i+2) & \text{ if $i \in \{1, 2, 3, 4, 5\}$, and}\cr
A_1 & \text{ if $i = s_{1'}$.}\cr
\endcases
$$  To show that $\psi_2$ is a homomorphism, it is enough to check the
defining relations of $W(D_6)$ are respected.  This is trivial, except as
regards the relations involving $s_{1'}$.  In the latter case, we need
to check the following identities: $$\eqalign{
A_1^2 &= I, \cr
(A_1 M^+_7(i+1, i+2))^2 &= I \text{ for $i \ne 2$},\cr
(A_1 M^+_7(3,4))^3 &= I.\cr
}$$  All of these are routine calculations, which can be done by hand.
\qed\enddemo

\proclaim{Lemma 4.7}
The subgroup $G_2$ of $GL_7({\Bbb R})$ that is generated by the matrices $$
M^+_7(2, 3), 
M^+_7(3, 4), 
M^+_7(4, 5), 
M^+_7(5, 6), 
M^+_7(6, 7)
$$ and $$
A_2 = \left( \matrix
0 & 1 & 1 & 0 & 0 & 0 & 0\cr
1/2 & 1/2 & -1/2 & 0 & 0 & 0 & 0\cr
1/2 & -1/2 & 1/2 & 0 & 0 & 0 & 0\cr
-1/2 & 1/2 & 1/2 & 1 & 0 & 0 & 0\cr
-1/2 & 1/2 & 1/2 & 0 & 1 & 0 & 0\cr
-1/2 & 1/2 & 1/2 & 0 & 0 & 1 & 0\cr
-1/2 & 1/2 & 1/2 & 0 & 0 & 0 & 1\cr
\endmatrix \right)
$$ is isomorphic to $W(D_6)$.  Furthermore, the matrices can be respectively
identified with the generators $$
s_1, s_2, s_3, s_4, s_5, s_{1'}
$$ of $W(D_6)$.
\endproclaim

\demo{Proof}
Regard the matrices in the statement as linear transformations from
${\Bbb R}^7$ to ${\Bbb R}^7$ with respect to the usual basis 
$e_1, e_2, \ldots, e_7$.  Define $$
v_1 = (1, 1/2, 1/2, 1/2, 1/2, 1/2, 1/2)^T
;$$ note that $v_1$ is an eigenvector for $A_2$ with eigenvalue $1$.
Clearly, the set $$
\{v_1, e_2, e_3, e_4, e_5, e_6, e_7\}
$$ is a basis for ${\Bbb R}^7$.  The matrices of the aforementioned linear
transformations with respect to the new basis are precisely the matrices
of Lemma 4.6, and the conclusion follows.
\qed\enddemo

\proclaim{Lemma 4.8}
The subgroup $G_3$ of $GL_7({\Bbb R})$ that is generated by the matrices $$
M^+_7(2, 3), 
M^+_7(3, 4), 
M^+_7(4, 5), 
M^+_7(5, 6), 
M^+_7(6, 7)
$$ and $$
A_3 = \left( \matrix
0 & 0 & 1 & 1 & 0 & 0 & 0\cr
-1/2 & 1 & 1/2 & 1/2 & 0 & 0 & 0\cr
1/2 & 0 & 1/2 & -1/2 & 0 & 0 & 0\cr
1/2 & 0 & -1/2 & 1/2 & 0 & 0 & 0\cr
-1/2 & 0 & 1/2 & 1/2 & 1 & 0 & 0\cr
-1/2 & 0 & 1/2 & 1/2 & 0 & 1 & 0\cr
-1/2 & 0 & 1/2 & 1/2 & 0 & 0 & 1\cr
\endmatrix \right)
$$ is isomorphic to $W(D_6)$.  Furthermore, the matrices can be respectively
identified with the elements $$
s_1, s_2, s_3, s_4, s_5, s_1 s_2 s_{1'} s_2 s_1
$$ of $W(D_6)$.
\endproclaim

\demo{Proof}
This follows from Lemma 4.7 and the identity $$
A_3 = M^+_7(2, 3) M^+_7(3, 4) A_2 M^+_7(3, 4) M^+_7(2, 3)
;$$ in particular, the groups $G_2$ and $G_3$ are equal.
\qed\enddemo

We now relate the above discussions of Coxeter groups back to our group $M_K$ of Definition 4.1.

\proclaim{Proposition 4.9}  The group $M_K$ is isomorphic to $W(D_6)$.  Under this isomorphism, the generators  $(23)_S$, $(34)_S$, $(45)_S$, $(56)_S$, $(67)_S$, and $(12)$ of $M_K$ correspond respectively to the elements $$
s_1, s_2, s_3, s_4, s_5, s_1 s_2 s_{1'} s_2 s_1
$$ of $W(D_6)$.
\endproclaim

\demo{Proof}  By the definitions of the relevant quantities, we have $S  M^+_7(i, i+1) S^{-1}=(i,i+1)_S$ for  $2\le i\le 6$, with $S$ and $(i,i+1)_S$ as in Proposition 3.3 and $M^+_7(i, i+1)$ as in Lemma 4.8.   Also, by direct computation, we find that $(12)=S A_3 S^{-1}$, with $A_3$ also as in Lemma 4.8.  But then $M_K=\{ S\tau S^{-1}\colon \tau\in G_3\}$, with $G_3$ as in that same lemma.  Since $G_3$ is isomorphic to $W(D_6)$, so is $M_K$, and we are done.  \qed  \enddemo

The final order of business, for this section, is to use the above results concerning $M_K$ to produce explicitly a complete set of coset representatives for $G_K\backslash M_K$.

\proclaim{Proposition 4.10}  A complete set of coset representatives for the   coset space $G_K\backslash M_K$ is given by$$C_K=\{p_j\colon 0\le j\le 15\}\cup \{n_j\colon 0\le j\le 15\}, $$where$$  p_{4q+r}= [(25)(36)(47)]_S(1234)^q[(25)(36)(47)]_S (1234)^r$$for $0\le q,r\le 3$, and $$ n_j=[(12)(234567)_S]^5 p_j$$for $0\le j\le 15$.\endproclaim

\demo{Proof} One checks directly (see Remark 4.11 below) that, if $\mu ,\nu\in C_K$ are unequal, then the first coordinate of $\mu\vec{x}$ is different from the first coordinate of $\nu\vec{x}$.   But by Remark 3.4, then, $\mu$ and $\nu$ must be inequivalent modulo $G_K$. 

The fact that the sixteen $p_j$'s, together with the sixteen $n_j$'s, provide a {\it complete} set of coset representatives now follows from the fact that $|M_K|/|G_K|=23040/720=32$.\qed \enddemo

\remark{Remark 4.11}  The transformation $$\iota=[(12)(234567)_S]^5$$figuring in the above definition of the $n_j$'s is an involution, and has the following effect:$$\iota\vec{x}=\biggl({1-a,1-b,1-c,1-d,\atop 2-e,2-f,2-g}\biggr)^T.$$This transformation $\iota$  corresponds, under the isomorphism described in Proposition 4.9, to the central involution of $W(D_6)$, which is the unique nontrivial element of the center of $W(D_6)$.

We also remark that, since  $[(25)(36)(47)]_S\in G_K$, we have $G_K [(25)(36)(47)]_S \mu= G_K\mu$ for any $\mu \in M_K$; therefore, the above proposition remains true if we omit  the leftmost ``$[(25)(36)(47)]_S$'' in the above definition of the $p_k$'s.  We include it for aesthetic reasons. Specifically, we do so to ensure that  (a)    $p_0$ is the identity matrix; and (b) the $p_j\vec{x}$'s and  $n_j\vec{x}$'s take particularly symmetric forms.  

Indeed,  consider the transformation $R$ of $\C^4$ defined by  $R(x,y,z,t)=(y,z,t,x)$.  Denote the images of the vectors $(a,b,c,d)$ and $(1,e,f,g)$ under the $j$th power transformation $R^j$ by  $(a_j,b_j,c_j,d_j)$ and $(1_j,e_j,f_j,g_j)$ respectively.  Then the above definitions of the $p_j$'s and $n_j$'s imply
$$  
{p_{4q+r}}\vec{x}=\biggl({1+a_r-1_q,1+b_r-1_q,1+c_r-1_q,1+d_r-1_q,\atop 1+e_q-1_q,1+f_q-1_q,1+g_q-1_q} \biggr)^T
 $$and 
$$  
{n_{4q+r}}\vec{x}=\biggl({1_q-a_r ,1_q-b_r,1_q-c_r,1_q-d_r,\atop 1+1_q-e_q,1+1_q-f_q,1+1_q-g_q} \biggr)^T
 $$for $0\le q,r\le 3$.  For example, $$p_{11}\vec{x}= p_{4\cdot 2+3}\vec{x}=\biggl({1+d-f,1+a-f,1+b-f,1+c-f,\atop 1+g-f,2-f,1+e-f}\biggr)^T$$ and $$n_{11}\vec{x}= n_{4\cdot 2+3}\vec{x}=\biggl({f-d,f-a,f-b,f-c,\atop 1+f-g,f,1+f-e}\biggr)^T.$$

The above explicit characterizations of our coset representatives explain our choice of notation:  in the $p_j$'s, the numerator parameters $a,b,c,d$ appear in a {\it positive} sense; in the $n_j$'s, they appear in a {\it negative} sense.
\endremark

We conclude this section with a labelling of the elements of $C_K$.

\definition{Definition 4.12}  Given $\mu\in C_K$, let $s(\mu)$ be the pre-image of $\mu$ under the isomorphism of Proposition 4.9.  Also, let $\Omega$ be as in Lemma 4.4 and $\omega_0$ as in Lemma 4.5. We define the {\it label attached to $\mu$} to be the six-digit binary number $b(\mu)$  (with evenly many ones) whose $k$th binary digit is $0$ if the $k$th coordinate of $s(\mu)\omega_0$ is
positive, and $1$ otherwise.   \enddefinition

For example, consider $p_1=(1234)=(12)(23)(34)=(12)(23)_S(34)_S\in C_K$.  By Proposition 4.9,  $s(p_1)=(s_1 s_2 s_{1'} s_2 s_1)(s_1)(s_2)=s_1 s_2 s_{1'} \in W(D_6)$.  But $$\displaylines{ s_1 s_2 s_{1'}w_0=M^+(1,2)  
M^+(2,3)
M^-(1, 2)(1,1,1,1,1,1)^T\cr=M^+(1,2)  
M^+(2,3)
 (-1,-1,1,1,1,1)^T=M^+(1,2)  
 (-1,1,-1,1,1,1)^T\cr= (1,-1,-1,1,1,1)^T,}$$so $b({p_1})=011000$.  

It also follows from the definition of the reversal $\iota$ that, for $0\le k\le 15$, the binary digits of $p_k$ and $n_k$ sum to $111111$.

 It is straightforward to compute the label of each element of $C_K$; we find: 

\settabs 4\columns
 
\+$b({p_0})\hskip4pt = 000000$,&
$b({p_1})\hskip4pt =011000$,&
$b({p_2})\hskip4pt =101000,$&
$b({p_3})\hskip4pt =110000$,\cr
\+$b({p_4})\hskip4pt=  000011$,&
$b({p_5})\hskip4pt =011011$,&
$b({p_6})\hskip4pt =101011,$&
$b({p_7})\hskip4pt =110011$,\cr
\+$b({p_8})\hskip4pt = 000101$,&
$b({p_9})\hskip4pt =011101$,&
$b({p_{10}})=101101,$&
$b({p_{11}})=110101$,\cr
\+$b({p_{12}})=  000110$,&
$b({p_{13}})=011110$,&
$b({p_{14}})=101110,$&
$b({p_{15}})=110110$;\cr
 \smallskip
\+$b({n_0})\hskip4pt = 111111$,&
$b({n_1})\hskip4pt =100111$,&
$b({n_2})\hskip4pt =010111,$&
$b({n_3})\hskip4pt =001111$,\cr
\+$b({n_4})\hskip4pt=  111100$,&
$b({n_5})\hskip4pt =100100$,&
$b({n_6})\hskip4pt =010100,$&
$b({n_7})\hskip4pt =001100$,\cr
\+$b({n_8})\hskip4pt = 111010$,&
$b({n_9})\hskip4pt =100010$,&
$b({n_{10}})=010010,$&
$b({n_{11}})=001010$,\cr
\+$b({n_{12}})=  111001$,&
$b({n_{13}})=100001$,&
$b({n_{14}})=010001,$&
$b({n_{15}})=001001$.\cr

\head 5. A metric  on $G_K\backslash M_K$ \endhead

\definition{Definition 5.1}
Let $\Omega$ be as in Lemma 4.4 (iii); let $\omega_1, \omega_2 \in \Omega$.  We define the {\it Hamming distance},
$d(\omega_1, \omega_2)$, between $\omega_1$ and $\omega_2$ to be the number
of coordinates at which the vectors $\omega_1$ and $\omega_2$ disagree.
\enddefinition

\remark{Remark 5.2}
Note that any two elements of $\Omega$ disagree at an even number (0, 2, 4
or 6) of coordinates.
\endremark

\proclaim{Lemma 5.3}
\item{}{\rm (i)} {The set $\Omega$ is a metric space with respect to Hamming
distance.}
\item{}{\rm (ii)} 
{For any $\omega_1, \omega_2 \in \Omega$, we have $\omega_1 . \omega_2 = 
6 - 2 d(\omega_1, \omega_2).$}
\item{}{\rm (iii)} {The action of $W(D_6)$ on $\Omega$ is by isometries with
respect to Hamming distance.}
\endproclaim

\demo{Proof}
Part (i) is a routine exercise using the definitions.  Part (ii) is also
routine, and may be checked case by case using Remark 5.2.

To prove part (iii), we note that the subgroup $\phi(W(D_6))$ (see
Lemma 4.4) of $GL_6({\Bbb R})$ consists of orthogonal matrices, so the
action of $W(D_6)$ respects scalar product.  By part (ii), the action also
respects Hamming distance.
\qed\enddemo

\definition{Definition 5.4}
Let $\Omega^{(3)}$ be the subset of the power set of $\Omega$ consisting of
all unordered triples $\{a, b, c\}$ of distinct elements of
$\Omega$ (i.e., $a \ne b \ne c \ne a$).
\enddefinition

\proclaim{Proposition 5.5}
\item{}{\rm (i)} {The group $W(D_6)$ acts on $\Omega^{(3)}$ diagonally via 
$g . \{a, b, c\} = \{g.a, g.b, g.c\}.$}
\item{}{\rm (ii)} {If $\{a_1, a_2, a_3\}$ and $\{b_1, b_2, b_3\}$ are two
elements of $\Omega^{(3)}$ for which $d(a_i, a_j) = d(b_i, b_j)$ for all
$i, j \in \{1, 2, 3\}$, then there exists $w \in W(D_6)$ such that 
$w(a_i) = b_i$ for all $i \in \{1, 2, 3\}$.}
\item{}{\rm (iii)} {Two elements $\{a, b, c\}$ and $\{a', b', c'\}$ of 
$\Omega^{(3)}$ are in the same $W(D_6)$-orbit if and only if the unordered
multisets $$\{d(a, b), d(b, c), d(a, c)\}$$ and 
$$\{d(a', b'), d(b', c'), d(a', c')\}$$ are equal.}
\item{}{\rm (iv)} {The group $W(D_6)$ has precisely $5$ orbits in its action
on $\Omega^{(3)}$.  These correspond to the multisets $\{2, 2, 2\}$, 
$\{4, 2, 2\}$, $\{4, 4, 2\}$, $\{4, 4, 4\}$ and $\{6, 4, 2\}$.}
\endproclaim

\demo{Proof}
Part (i) is an easy exercise.

Assume the hypotheses of (ii).  Since $W$ acts transitively on $\Omega$ by 
Lemma 4.4 (iii), we may assume without loss of generality that
$a_1 = b_1= \omega_0$, where $\omega_0$ is as in Lemma
4.5 (ii).  Let $T_2$ (respectively, $T'_2, T_3, T'_3$) be the subset of 
$T = \{1, 2, 3, 4, 5, 6\}$ corresponding to the occurrences of $-1$ in $a_2$
(respectively, $b_2, a_3, b_3$).  The hypotheses imply that $$
|T_2| = |T'_2| = d(a_1, a_2) = d(b_1, b_2)
,$$ $$
|T_3| = |T'_3| = d(a_1, a_3) = d(b_1, b_3)
$$ and $$
|T_2 \ \Delta \ T_3| =
|T'_2 \ \Delta \ T'_3| =
d(a_2, a_3) = d(b_2, b_3)
,$$ where $\Delta$ is the symmetric difference operator.  It follows that
there is a permutation $w$ of $T$ with the properties that 
(a) $w(T_2) = T'_2$, (b) $w(T_3) = T'_3$ and (c) $w(T_2 \cap T_3) = 
T'_2 \cap T'_3$.  We may regard this 
permutation as an element of $W(D_6)$ fixing $\omega_0$ by Lemma
4.4 (ii), so that $w(a_1) = b_1$, and we also have $w(a_2) = b_2$ and 
$w(a_3) = b_3$ as required, completing the proof of (ii).

The ``if'' direction of (iii) follows from (ii), and the ``only if''
direction follows from the fact, proved in Lemma 5.3 (iii), that $W$
acts on $\Omega$ by isometries.

To prove (iv), it is a routine exercise to construct triples in $\Omega^{(3)}$
corresponding to any of the desired multisets; for example, the triple $$
(+1, +1, +1, +1, +1, +1)^T,
(-1, -1, +1, +1, +1, +1)^T,
(+1, +1, +1, +1, -1, -1)^T,
$$ corresponds to the multiset $\{4, 2, 2\}$.  Conversely, we need to show
that the five multisets listed are exhaustive; the conclusion will then
follow from (iii).  If $a, b, c \in \Omega$ satisfy $d(a, b) = 6$ and
$c \not\in \{a, b\}$, then it must be the case that $c$ disagrees with one
of $a$ and $b$ in two places and disagrees with the other in four places.
It follows that the only multiset containing $6$ is $\{6, 4, 2\}$.  The
only other possible distances between two distinct 
points of $\Omega$ are $2$ and $4$,
and they occur in all possible combinations, completing the proof.
\qed\enddemo

\remark{Remark 5.6}
Proposition 5.5 does more than describe the orbits of $W(D_6)$ on
$\Omega^{(3)}$.  For example, if $\{a, b, c\} \in \Omega^{(3)}$ has the
property that $d(a, b) = d(a, c) = 4$ and $d(b, c) = 2$, then Proposition
5.5 (ii) shows that there exists $w \in W(D_6)$ such that $w(a) = a$,
$w(b) = c$ and $w(c) = b$.
\endremark

The Hamming distance described above, together with the labeling of cosets defined in the previous section, provide a notion of Hamming distance on our coset space $G_K\backslash M_K$.  Namely: the Hamming distance between two cosets, represented by elements $\mu,\nu\in C_K$ respectively, is just the number of digits at which the labels for $\mu$ and $\nu$ disagree.  

For example, the cosets represented by $p_3$ and  $n_{12}$ have labels   110000 and 111001 respectively; the Hamming distance between them is therefore equal to $$|1-1|+|1-1|+|0-1|+|0-0|+|0-0|+|0-1|=2.$$

Possible values for the Hamming distance are $0$, $2$, $4$, and $6$, cf. Remark 5.2 above.  The Hamming distance between $\mu,\nu,\in C_K$ will be 6 if and
only if $\mu=\iota \nu$, with $\iota$ the reversal involution described in Remark 4.11 above.

 \head \S6.  Three-term relations \endhead

The Hamming distance on $G_K\backslash M_K$, defined in the previous section, provides a convenient classification of triples $(\mu_1,\mu_2,\mu_3)$ of elements of $C_K$.  This classification will be central to our study, below, of three-term relations for $K(a;b,c,d;e,f,g)$. 

\definition{Definition 6.1} The Hamming type of a triple $(\mu_1,\mu_2,\mu_3)$ of elements of $C_K$ is  defined to be the three-digit integer $abc$, where $a$ is the shortest of the Hamming distances among $\mu_1,\mu_2$ and $\mu_3$; $b$ is the next shortest; and $c$ is the longest.  \enddefinition 

Thus, the possible Hamming types are $222$, $224$, $244$, $444$, and $246$.  (See the proof of Proposition 5.5.)  For example:  we have noted previously that $p_3$ and  $n_{12}$ are separated by a Hamming distance of 2.  The first of these cosets has Hamming distance 2 from $p_0$, and the second has Hamming distance $4$.  So the triple $(p_0,p_3,n_{12})$ is of Hamming type 224.

Also, for $\mu\in C_K$ and $V$ the affine hyperplane of (2.1), we define functions $K_\mu\colon V \to\C$ by $K_\mu(\vec{x})=K(\mu\vec{x})$.  We call a relation among $K_{\mu_1},K_{\mu_2}$, and $K_{\mu_3}$ an ``$abc$ relation'' if the triple $(\mu_1,\mu_2,\mu_3)$ is of Hamming type $abc$.

We are now ready to state, and prove, our main results concerning three-term relations for $K(a;b,c,d;e,f,g)$.  
It is also convenient to introduce the following functions, which
occur frequently in our identities.

\definition{Definition 6.2}
Suppose that $\vec{x}=(a,b,c,d,e,f,g)^T\in\ V$.
We define $$\aligned&\alpha(\vec{x})=\frac{\sin\pi(c-b)}{\gg{e-a}\gg{f-a}\gg{g-a}},\\&\beta(\vec{x}) =\frac{1}{\gg{a}\gg{1+b-e}\gg{f-c} \gg{f-d}\gg{g-c}\gg{g-d}},\\&\gamma(\vec{x}) = \frac{\sin\pi(e-b)}{\sin\pi(c-b)}\biggl[{\sin\pi(a-c)\sin\pi(e-a-b)\sin\pi(f-b)\sin\pi(g-b)\atop-
\sin\pi(a-b)\sin\pi(e-a-c)\sin\pi(f-c)\sin\pi(g-c) }\biggr] .\endaligned$$
\enddefinition

\proclaim{Proposition 6.3}
We have the 222 relation$$\alpha(\vec{x})K_{p_0}(\vec{x})
+\alpha((123)\vec{x})K_{p_1}(\vec{x})+\alpha((132)\vec{x})K_{p_2}(\vec{x})=0,\eqno{(6.1)}$$or in other words,$$\sum_{j=0}^2 \alpha((123)^j\vec{x}) K ((123)^j\vec{x})=0.$$
\endproclaim

\demo{Proof}
By Proposition 3.3, $K(\vec{x})$ is invariant under the substitution $\vec{x}\to (a,g-c,g-b,d,1+a+d-f,1+a+d-e,g)^T.$  So we may apply this substitution to the right-hand side of Proposition 3.1; the result is
$$\displaylines{  K(\vec{x})
 =\frac{1}{ {
\gg{f-d}\gg{e-d}\gg{g-c}\gg{g-b}\gg{e-b}\gg{e-c}\gg{f-b}\gg{f-c}
\gg{g-d} }}
\cr\cdot{1\over 2\pi i}
\int {\gg{g-d+t}\gg{f-d+t}\gg{e-d+t}\gg{d-c-t}
\gg{d-b-t}\gg{-t}\over\gg{d-t}\gg{1+a-d+t}} \,dt 
\cr=\frac{1}{ {
\gg{e-b}\gg{e-c}\gg{e-d}\gg{f-b}\gg{f-c}\gg{f-d}
\gg{g-b}\gg{g-c}\gg{g-d} }}
\cr\cdot{1\over 2\pi i}
\int {\gg{g+t}\gg{f+t}\gg{e+t}\gg{-d-t}\gg{-c-t}
\gg{-b-t}\over\gg{-t}\gg{1+a+t}} \,dt  ,}$$the last step by some rearranging and the substitution $t\to t+d$.  But then  $$\displaylines{ \sum_{j=0}^2 \alpha((123)^j\vec{x}) K ((123)^j\vec{x}) 
\cr=\frac{1}{\biggl[{\ds\gg{e-a}\gg{f-a}\gg{g-a}\gg{e-b}\gg{f-b}\gg{g-b}\atop\ds\cdot \gg{e-c}\gg{f-c}\gg{g-c}\gg{e-d}\gg{f-d}\gg{g-d}}\biggr]}
\cr\cdot{1\over 2\pi i}
\int {\gg{g+t}\gg{f+t}\gg{e+t} \gg{-d-t}\over\gg{-t}} \biggl[\frac{\sin\pi (c-b)\gg{-c-t}
\gg{-b-t} }{\gg{1+a+t}}\cr +\frac{\sin\pi (a-c)\gg{-a-t}
\gg{-c-t} }{\gg{1+b+t}}+\frac{\sin\pi (b-a)\gg{-b-t}
\gg{-a-t} }{\gg{1+c+t}}\biggr]\,dt
\cr}$$ $$\displaylines{=\frac{1}{\biggl[{\ds\gg{e-a}\gg{f-a}\gg{g-a}\gg{e-b}\gg{f-b}\gg{g-b}\atop\ds\cdot \gg{e-c}\gg{f-c}\gg{g-c}\gg{e-d}\gg{f-d}\gg{g-d}}\biggr]}
\cr\cdot{1\over 2\pi i}
\int {\gg{g+t}\gg{f+t}\gg{e+t}\gg{-d-t} \gg{-c-t}
\gg{-b-t}\gg{-a-t}\over\gg{-t}} \cr\cdot\frac{1}{\pi} \bigl[ \sin\pi (c-b)\sin\pi (-a-t)  + \sin\pi (a-c)\sin\pi(-b-t)  + \sin\pi (b-a)\sin\pi (-c-t)\bigr]\,dt,}$$the last step by (2.6).  The quantity in square-brackets on the right-hand side is seen to equal zero, by straightforward application of trigonometric identities, which completes the proof.
\qed\enddemo

\proclaim{Proposition 6.4}
We have the 224 relation$$ \pi^3 \sin\pi(e-a-b) \beta(\sigma \vec{x})K_{p_0}(\vec{x})+\gamma(\vec{x})K_{p_1}(\vec{x})-\pi^3 \sin\pi(a-b) \beta( \vec{x})K_{n_4}(\vec{x})=0,$$where$$ \sigma\vec{x}=\biggl({ 1+b-e, b,f-c,g-c,\atop 1+a+b-e,1+b-c,1+b+d-e}\biggr)^T   $$(note that $\sigma^3\vec{x}=\vec{x}$).
\endproclaim

\demo{Proof}
Into (6.1) we substitute $$\vec{x}\to 
\biggl({ b,e-a,c,e-d,\atop e,1+b+c-g,1+b+c-f}\biggr)^T .$$ It is readily checked that this transformation takes $K_{p_0}$ to
$K_{p_1}$; takes $K_{p_1}$ to $ K_{n_4}$; and takes $K_{p_2}$ to itself.  So (6.1) yields $$\aligned \frac{\sin\pi (a+c-e)\,K_{p_1}(\vec{x})}{\gg{e-b}\gg{1+c-g}\gg{1+c-f} }&+\frac{\sin\pi (b-c)\,K_{n_4}(\vec{x})}{\gg{a}\gg{f-d}\gg{g-d} }\\&+\frac{\sin\pi (e-a-b)\,K_{p_2}(\vec{x})}{\gg{e-c}\gg{1+b-g}\gg{1+b-f} }=0.\endaligned\eqno{(6.2)}$$

We  multiply  (6.1)   by $ \pi^3\sin\pi(e-a-b)/( \sin\pi(c-b) \gg{1+b-e}\gg{1+b-f} \gg{1+b-g}) $, multiply  (6.2) by 
$  \pi^3\sin\pi (b-a) /(\sin\pi(c-b) \gg{1+b-e} \gg{f-c}\gg{g-c} ) $, and subtract; applying (2.6) to the resulting coefficient of $K_{p_1}(\vec{x})$, we find that $$\aligned&\frac{\pi^3 \sin\pi(e-a-b)\,K_{p_0}(\vec{x})}{\gg{e-a}\gg{f-a}\gg{g-a}  \gg{1+b-e}\gg{1+b-f}\gg{1+b-g}}\\+&  \frac{\sin\pi(e-b)}{  \sin\pi(c-b)}\biggl[{\sin\pi(a-c)\sin\pi(e-a-b)\sin\pi(f-b)\sin\pi(g-b)\atop-
\sin\pi(a-b)\sin\pi(e-a-c)\sin\pi(f-c)\sin\pi(g-c) }\biggr]K_{p_1}(\vec{x})\\-&\frac{ \pi^3\sin\pi(a-b)\,K_{n_4}(\vec{x})}{\gg{a} \gg{1+b-e} \gg{f-c}\gg{g-c}\gg{f-d}\gg{g-d} }=0,\endaligned\eqno{(6.3)}$$ as required.
\qed\enddemo

\proclaim{Proposition 6.5}
We have the 244 relation$$\displaylines{\sin \pi a \sin\pi(f-e)\sin\pi(g-c)     \sin\pi  (g-d)  \beta(\sigma \vec{x})K_{p_0}(\vec{x})\cr +\gamma((56)\vec{x})\beta( \vec{x})K_{n_4}(\vec{x})-\gamma( \vec{x})\beta( (56)\vec{x})K_{n_8}(\vec{x}) =0.}$$
\endproclaim

\demo{Proof}
We interchange $e$ and $f$ in (6.3) to get a relation among $K_{p_0}$, $K_{p_1}$, and $K_{n_8}$.    From this relation and from (6.3) itself, we may now eliminate $K_{p_1}$.  In doing so, we note that two of the terms in the resulting coefficient of $K_{p_0}$ cancel, and that the remaining two terms may be combined using   the trigonometric identity
$$ \aligned &\sin\pi(f-b) \sin\pi(f-a-c)   \sin\pi(e-c)\sin\pi (e-a-b) 
\\-& \sin\pi(e-b) \sin\pi(e-a-c)   \sin\pi(f-c)\sin\pi (f-a-b)  \\=  &\sin \pi a     \sin\pi  (c-b)  \sin \pi(f-e) \sin\pi (g-d)\endaligned\eqno{(6.4)} $$
for $e+f+g-a-b-c-d=1$.  The final result, after some simplification, is $$\aligned  & \frac{ \sin\pi a\sin\pi(f-e) \sin\pi(g-c)\sin\pi(g-d)    }{ \gg{e-a}\gg{f-a}\gg{g-a} \gg{1+b-e}\gg{1+b-f} \gg{1+b-g} }\,K_{p_0}(\vec{x})
 \\+& \frac{\ds\frac{\sin\pi(f-b)}{\sin\pi(c-b)} \biggl[{\ds\sin\pi(a-c)\sin\pi(f-a-b)\sin\pi(e-b)\sin\pi(g-b)\atop\ds-
\sin\pi(a-b)\sin\pi(f-a-c)\sin\pi(e-c)\sin\pi(g-c) }\biggr] }{  \gg{a} \gg{1+b-e}\gg{f-c} \gg{f-d}   \gg{g-c} \gg{g-d}}\, K_{n_4}(\vec{x})
 \\-&  \frac{\ds\frac{\sin\pi(e-b)}{\sin\pi(c-b)}  \biggl[{\ds\sin\pi(a-c)\sin\pi(e-a-b)\sin\pi(f-b)\sin\pi(g-b)\atop\ds-
\sin\pi(a-b)\sin\pi(e-a-c)\sin\pi(f-c)\sin\pi(g-c) }\biggr] }
{ \gg{a} \gg{1+b-f} \gg{e-c} \gg{e-d}  \gg{g-c} \gg{g-d}} \,K_{n_8}(\vec{x})=0,\endaligned   $$as required.
\qed\enddemo

\proclaim{Proposition 6.6}
We have the 444 relation
$$\displaylines{  \gamma(\sigma^{2}\vec{x})\beta(\sigma\vec{x})K_{p_0}(\vec{x})+\gamma(\sigma\vec{x}) \beta( \vec{x})K_{n_4}(\vec{x})+ \gamma(\vec{x})\beta( \sigma^{-1}\vec{x})K_{p_5}(\vec{x}) =0,}$$or in other words,$$\sum_{j=0}^2\gamma(\sigma^{2-j}\vec{x}) \beta(\sigma^{1-j}\vec{x})K(\sigma^{-j}\vec{x})=0.$$
\endproclaim

\demo{Proof}
We apply, to (6.3), the substitution    $$\vec{x}\to \biggl({a,b,f-c,f-d,\atop 1+a+b-e,f,1+a+b-g}\biggr)^T. $$ This substitution takes  $K_{p_0}$ to itself, takes   $K_{p_1}$ to itself, and takes  $K_{n_4}$ to $K_{p_5}$.  So we obtain a relation among $K_{p_0}$, $K_{p_1}$, and $K_{p_5}$; from this relation and from (6.3) itself  we eliminate $K_{p_1}$, to obtain a relation among $K_{p_0}$, $K_{n_4}$, and $K_{p_5}$.

The calculations here are much as they were in the proof of Proposition 6.5.  
The salient difference in the present case is that the coefficient of $K_{p_0}$ reduces, by way of some trigonometric identities like (6.4), to a {\it binomial} instead of a monomial in gamma and sine functions.  The upshot (which again makes use of (2.6)) is the following 444 relation: 
$$\aligned &\frac{ \ds\frac{\sin\pi a}{\sin \pi f } \ds \biggl[{\sin \pi (c+d-g) \sin \pi e \,\sin  \pi (f-c) \sin \pi (f-d) \atop-\sin  \pi (e-a-b) \sin \pi (f-e) \sin \pi c \,\sin\pi d}\biggr] }{  \gg{e-a}   \gg{f-a}\gg{g-a}  \gg{1+b-e}\gg{1+b-f}\gg{1+b-g}  }\, K_{p_0}(\vec{x})
 \\+&
   \frac{\ds\frac{\sin\pi(1+a-e)}{\sin\pi (f-b-c)}\ds\biggl[ {
\sin\pi (a-b)\sin\pi (g-a-d) \sin\pi c\,\sin\pi(e-d)\atop
- \sin\pi (a+c-f)\sin\pi e \,\sin\pi(f-b)\sin\pi(g-a) } \biggr] }{ \gg{a}\gg{1+b-e}  \gg{f-c}\gg{f-d}\gg{g-c} \gg{g-d}  } 
 \,  K_{n_4}(\vec{x})
\\+& \frac{ \ds\frac{\sin\pi(e-b)}{\sin\pi (c-b)}\ds\biggl[ {  \sin\pi (a-c)\sin\pi (e-a-b)\sin\pi(f-b)\sin\pi(g-b)\atop -
\sin\pi (a-b)\sin\pi (e-a-c) \sin\pi(f-c)\sin\pi(g-c) } \biggr]  }{ \gg{a} \gg{c}\gg{d}  \gg{e-a} \gg{e-c}\gg{e-d}}
\,  
K_{p_5}(\vec{x}) =0,\endaligned $$as required.
\qed\enddemo

\proclaim{Proposition 6.7}
We have the 246 relation$$\aligned\frac{\sin\pi(e-b)}{\sin\pi(c-b)} &\biggl[{\sin\pi(f-c)\sin\pi(g-c)\sin\pi(e-a-c)\gamma(\sigma
(12)(34) \vec{x})\atop-\sin\pi(f-b)\sin\pi(g-b)\sin\pi(e-a-b)\gamma(\sigma(1342) \vec{x})}\biggr]K_{p_0}(\vec{x})\\+&{\pi^6\sin\pi e}\, \beta(\vec{x})\beta(\sigma(12)\vec{x})K_{n_4}(\vec{x})
 + \pi^3 \gamma(\vec{x})\beta(\sigma^{-1}(12)\vec{x})K_{p_4}(\vec{x})=0.\endaligned$$
\endproclaim

\demo{Proof}
Into (6.3) we substitute  $$\vec{x}\to \biggl({b,a,f-c,f-d,\atop 1+a+b-e,f,1+a+b-g}\biggr)^T .$$ This substitution takes  $K_{p_0}$ to $K_{p_1}$, takes $K_{p_1}$ to $K_{p_0}$, and takes $K_{n_4}$ to $K_{p_4}$.  So we obtain a relation among $K_{p_0}$, $K_{p_1}$, and $K_{p_4}$; from this relation and from (6.3) itself, we eliminate $K_{p_1}$ to obtain a relation among $K_{p_0}$, $K_{n_4}$, and $K_{p_4}$, which is of type 246.

The above procedure yields a relation in which the coefficient of $K_{p_0}$ contains, {\it a priori}, five monomials in gamma and sine functions.  But a trigonometric identity similar to (6.4) allows us to combine two of these monomials into one. 
The final result is
$$\aligned & 
\frac{\sin\pi(e-b)}{\sin \pi (c-b)   }\biggl(
\sin \pi  (f-b)\sin \pi (g-b) \sin \pi(e-a-b)  
\\\cdot& \frac{ \sin   \pi (e-c ) }{ \sin \pi (f-a-d) }\biggl[{  \sin\pi (c+d-f) \sin \pi  e \,\sin
   \pi (f-a) \sin \pi (g-c) 
\atop  -\sin \pi  (c-a)  \sin  \pi (g-b-c)   \sin \pi d      \,\sin \pi (e-b) 
  } \biggr]
\\-& \sin \pi (f-c)  \sin \pi (g-c)   \sin  \pi (e-a-c)
\\\cdot&\frac{ \sin   \pi (e-b)}{\sin \pi (f-a-d) }\biggl[{  \sin\pi (b+d-f) \sin \pi  e \,\sin
   \pi (f-a) \sin \pi (g-b) 
\atop  -\sin \pi  (b-a)  \sin  \pi (g-b-c)   \sin \pi d      \,\sin \pi (e-c) 
  } \biggr]\biggr)K_{p_0}(\vec{x})\\+&
 \frac{\pi^6\sin\pi e  }{ \biggl[{\displaystyle \gg{1+a-e} \gg{1+a-f}\gg{1+a-g} \gg{e-b}\gg{f-b}\gg{g-b}
\atop\cdot\displaystyle\gg{a}\gg{1+b-e} \gg{f-c}\gg{g-c} \gg{f-d}\gg{g-d}}\biggr] }  
 \, K_{n_4}(\vec{x}) \\+&
  \frac{\pi^3 }{  \gg{b} \gg{c}\gg{d}  \gg{e-b}\gg{e-c}\gg{e-d} } 
\\\cdot &\displaystyle\frac{ \sin\pi(e-b)}{\sin\pi(c-b)} \biggl[ { \sin\pi (a-c)\sin\pi (e-a-b) \sin\pi(f-b)\sin\pi(g-b)\atop -
\sin\pi (a-b)\sin\pi (e-a-c) \sin\pi(f-c)\sin\pi(g-c) } \biggr]
 K_{p_4}(\vec{x}) =0.\endaligned $$ This completes the proof.
\qed\enddemo

\proclaim{Theorem 6.8}

{\rm (i)}  Let $\mu_1,\mu_2,\mu_3$ be elements of $M_K$, such that no two of the $\mu_j$'s are in the same right coset of $G_K$ in $M_K$.  Then there is a relation of the form
$$ \gamma_1  K_{\mu_1} + \gamma_2 \,K_{\mu_2} + \gamma_3  K_{\mu_3}  =0,\eqno{(6.5)} $$where  $\gamma_1$, $\gamma_2$, $\gamma_3$ are entire, rational combinations of gamma and sine functions, whose arguments are $\Z-$linear combinations of $a,b,c,d,e,f,g$.

{\rm (ii)}  For any $\ell\in\{1,2,3\}$, write $\{1,2,3\}-\{\ell\}=\{j,k\}$.  Then in a relation of the form {\rm (6.5)}, each coefficient $\gamma_\ell$  may written as a sum of $2^{n-1}$ monomials in gamma and sine functions, where $2n$ is the Hamming distance between $\mu_j$ and $\mu_k$.

{\rm (iii)}  If $(\mu_1,\mu_2,\mu_3)$ and $(\nu_1,\nu_2,\nu_3)$ are triples of the same Hamming type, then a three-term relation among $K_{\mu_1},K_{\mu_2}$, and $K_{\mu_3}$ can be transformed into one involving $K_{\nu_1},K_{\nu_2}$, and $K_{\nu_3}$ by the application of a single change of variable$$\vec{x}\to \rho\vec{x}\quad(\rho\in M_K)$$to all elements (including the coefficients) of the first relation.
\endproclaim

\demo{Proof} The assertions follow by combining Proposition 5.5 with 
Propositions 6.3--6.7. \qed\enddemo

\leftheadtext{} \rightheadtext{}\vfill\eject

\end